\documentclass{article}

\usepackage[utf8]{inputenc}
\usepackage{authblk}
\usepackage{amssymb,latexsym,amsmath}
\usepackage{amsthm}
\usepackage{dsfont}
\usepackage{amsfonts}
\usepackage{mathrsfs}
\usepackage{multirow}
\usepackage[english]{babel}
\usepackage{array}
\usepackage{tikz-cd}
\usepackage{tikz-cd}
\usepackage{mathtools}
\usepackage{pst-node}
\usepackage{amscd}

\usepackage{tikz-cd}

\usepackage[verbose]{placeins}
\usepackage{graphicx}
\graphicspath{ {/admin/desktop/} }

\newcommand{\R}{\mathbb{R}}

\newcommand{\Z}{\mathbb{Z}}

\newtheorem{teorema}{Theorem}[section]
\newtheorem{lema}{Lemma}[section]

\DeclareMathOperator{\GL}{GL}
\title{An explicit section of the Laudenbach exact sequence of the mapping class group of    connect sums of $S^2 \times S^1$}
\author{Jorge Andres Robinson Arrieta }
\date{October ,18, 2023}

\begin{document}

\maketitle

\begin{abstract}\noindent
 Laudenbach proved that  the mapping class group of the connect sum of $n$ copies of $S^2 \times S^1$   is an extension of $Out(F_n)$ by a finite group. Brendle-Broaddus-Putman proved that this exact sequence splits. We provide an explicit section $s$ of this split exact sequence.
\end{abstract}

\section{Introduction } 
Let  $M_n$ be the connected sum of $n$ copies of $S^2 \times S^1$  equipped with a basepoint $x_0$. $Mod(M_n)$ is defined to be the group of isotopy classes of orientation-preserving diffeomorphisms of $M_n$. We fix an isomorphism $\pi_1(M_n,x_0)\cong F_n$, where $F_n$ is the free group of rank $n$.
In [1,2], Laudenbach proved that that there exists a short exact sequence
$$1 \xrightarrow{} Twist(M_n)\xrightarrow{} Mod(M_n)
\xrightarrow{\rho} Out(F_n) \xrightarrow{}1,$$ 
 where $Twist(M_n) \cong (\Z/2)^{n}$  is generated by the sphere twists about
the core spheres $S^{2} \times *$. Brendle-Broaddus-Putman proved
in [3]  that this short exact sequence splits. In particular, they construct a crossed homomorphism $\mathfrak{T}: Mod(M_n)\xrightarrow{} Twist(M_n)$ that restricts to
the identity on $Twist(M_n)$. This determines a section $s:Out(F_n)\xrightarrow{} Mod(M_n)$ of $\rho$ , given by $s([\phi])=\mathfrak{T}([f^{-1}])[f]$, where $f$ is a diffeomorphism of $M_n$ with $\rho([f])=[\phi]$. The purpose of this paper is to provide a formula for the section $s$ explicitly.  In order to do that, we compute $s$ for the Nielsen generators of $Out(F_n)$ given in [4]. We first describe explicit diffeomorphisms for each of the elements of the Nielsen generating set for $Out(F_n)$.  Our computation shows that $\mathfrak{T}$ is trivial for these lifts. Our main result is the following:

\begin{teorema}
    The map $s:Out(F_n)\xrightarrow{} Mod(M_n)$ that on the Nielsen generators $[R_{i, j}]$, and $[I_j]$, for $1\le i,j\le n$ and $i\neq j$, given by:
    $$ s([R_{i,j}])=[F_{i,j}], \text{ and } s([I_j])=[G_j],$$

    is a section of $\rho$, where $F_{i,j}$, and $G_j$ are diffeomorphisms of $M_n$ defined in the section below. 
\end{teorema}

\section{Construction of the maps $F_{i,j}$, and $G_j$ }
For $1\le i\le n$, choose loops $a_i$ based at $x_0$ that generate the fundamental group of  $M_n$. In [4, proposition 4.1], it is shown that $Out(F_n)$ is generated by the classes $[R_{i,j}],$ and $ [I_j]$, for $1\le i,j\le n$ and $i\neq j$, where:

\begin{equation*}
    R_{i,j}(a_k) = \begin{cases}
              a_ka_j & \text{if } k=i\\
               a_k & \text{if } k\neq i\\
               
          \end{cases}, \text{ and }
I_j(a_k) = \begin{cases}
              a_k^{-1} & \text{if } k=j\\
               a_k & \text{if } k\neq j.\\
               
          \end{cases}
\end{equation*}

We want to obtain diffeomorphisms $F_{i,j},$ and $ G_j$ of $M_n$ such that $\rho([F_{i,j}])=[R_{i,j}]$ and $\rho([G_j])= [I_j]$.

$M_n$ can be described by removing $2n$ open balls of $S^{3}$, and then gluing the boundary spheres of these balls in pairs. The resulting boundary spheres correspond to the core spheres of the $n$ summands of $S^{2}\times S^{1}$ in  $M_n$. Let $A_i$ denote the core sphere  of the $ith$ summand  $S^{2}\times S^{1}$ of $M_n$. Let $A_i^{-}$ and $A_i^{+}$ denote the two boundary spheres that were identified in $S^{3}$ minus $2n$ open balls to gave rise  the boundary sphere $A_i$ in $M_n$.  Define $a_i$ as the equivalence class of the curve starting at the base point that reaches $A_i^{-}$ and comes back trough $A_i^{+}$, and then reaches the base point, without intersecting the other boundary spheres. Then, $\{a_1,...,a_n\}$ forms a basis for $\pi_1(M_n,x_0)$.    Choose a subset $N_{i,j}$ of $M_n$ diffeomorphic to $D^2\times S^1$ minus an open ball with boundary  $A_i$, which is contained in $\{(x,y)|x^2+y^2<1/9\}\times S^1$, where $*\times S^1$ is freely homotopic to $a_j$, as depicted in  Figure 1 for the case $n=3, i=1$, and $j=2$. 

\begin{figure}[h]

\includegraphics[width=0.5\textwidth]{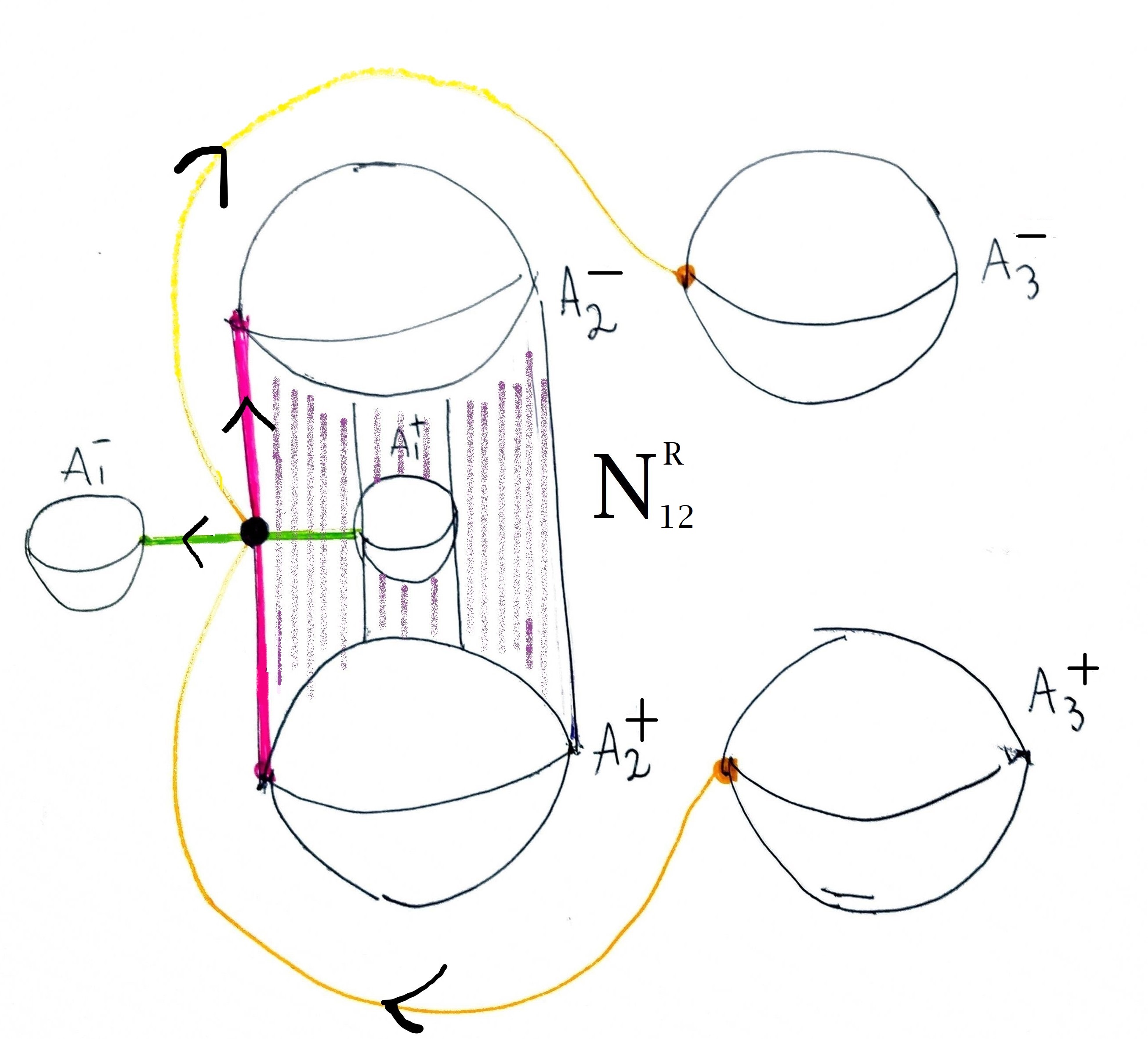}
\caption{$a_1, a_2$, and $a_3$ are depicted in green, red and yellow respectively. The neighborhood $N_{1,2}$ is depicted in purple.}
\end{figure}
\FloatBarrier

For the case of $I_j$, choose a subset $P_j$ of $M_n$ diffeomorphic to $B^3=\{(x,y,z)|x^2+y^2+z^2\le 1\}$ minus two open balls with boundaries  $A_j$, which are contained in $\{(x,y,z)|x^2+y^2+z^2<1/9\}$, and are symmetric about the $z-$axis. Parametrize $a_j$ in such a way that $a_j(t)\in \{(x,y,z)|x^2+y^2+z^2<1/9\}$ for $1/3\le t \le 2/3$, and $a_j(t)\in \{(x,y,z)|1/9<x^2+y^2+z^2<1\}$ for $0\le t \le 1/3$ and $2/3\le t\le 1$. we also homotope $a_j$ so that $a_j(t)$ lives in the $z-$axis for $t$ as in the last case.   Figure 2 depicts this for the case $j=1$. 

\begin{figure}[hbt!]

\includegraphics[width=0.6\textwidth]{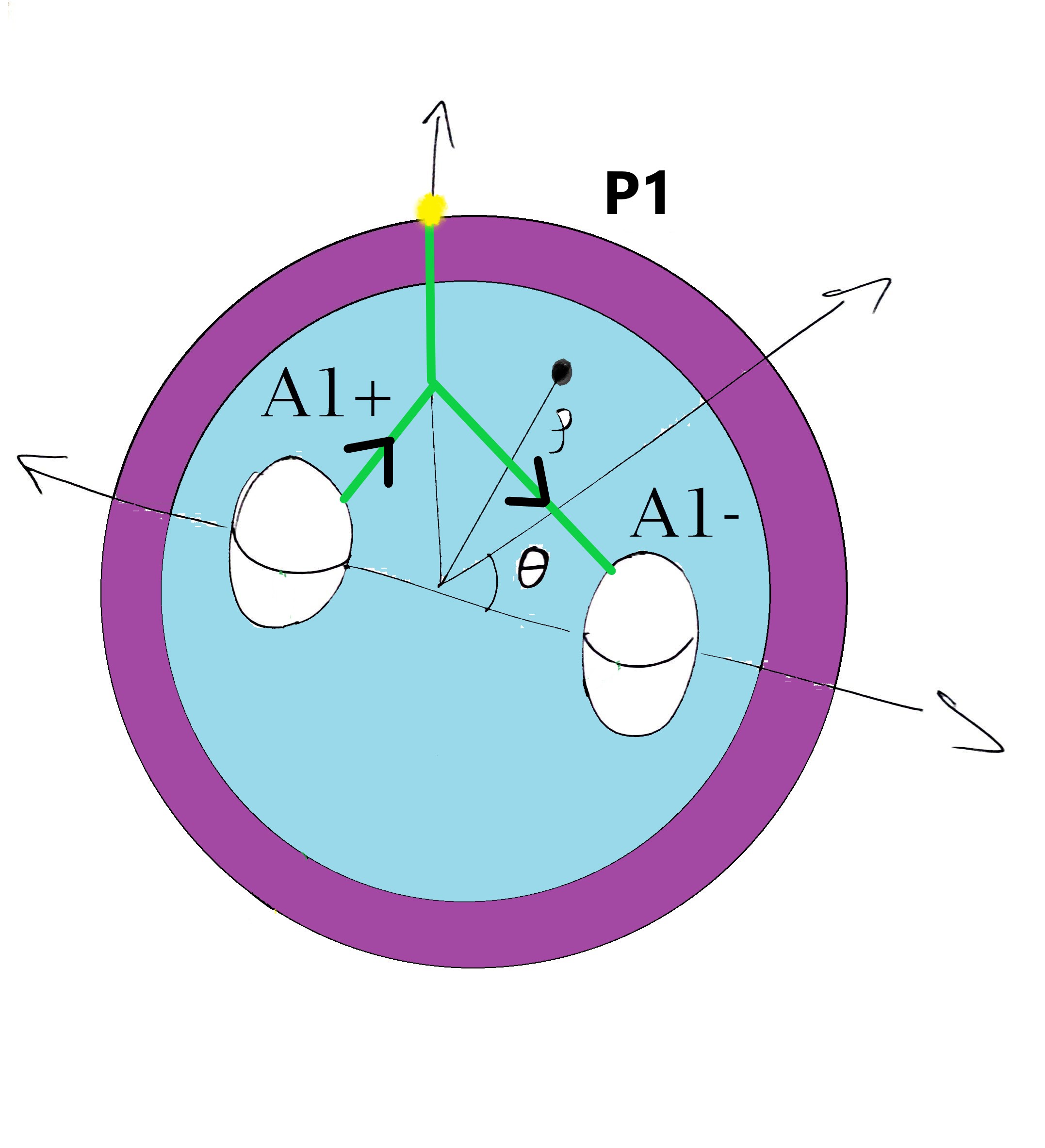}
\caption{$a_1$ is depicted in green}
\end{figure}
 \FloatBarrier

Construct a smooth  function $\psi:[0,1]\xrightarrow{} [0,1]$ with $\psi(r)=1$ on [0,1/3], $supp(\psi(r))\subseteq [0,2/3)$, and decreasing, so  that $\psi^{\prime}(r)\le 0$.

Define $f_{i,j}:N_{i,j}\xrightarrow{}N_{i,j}$ by
$$ f(x,y,e^{2\pi \textrm{i}\theta})=(x,y, e^{2\pi \textrm{i}\left[\theta+\psi\left(\sqrt{x^{2}+y^{2}}\right)\right]}).$$

Then $f_{i,j}$ is a diffeomorphism of $N_{i,j}$.

$F_{i,j}$ is defined  by :

\[   F_{i,j}(p)=\left\{
\begin{array}{ll}
       f_{i,j}(p) & \text{ if } p\in N_{i,j}  \\
     p & p\in M_n-N_{i,j}. \\
      
\end{array}
\right. \]

If $p$ has spherical coordinates $(\theta, \varphi,r)$, define $G_j$  by :

\[   G_j(\theta, \varphi,r)=\left\{
\begin{array}{ll}
       (\theta+\psi(r)\pi, \varphi, r) & \text{ if } p\in P_j  \\
     p & p\in M_n-P_j .\\
      
\end{array} 
\right. \]

To see that $F_{i,j}$ actually realizes $R_{i,j}$, consider what $F_{i,j}$ does to the $a_k's$, as depicted in figure 3 for the case $n=3,i=1,$ and $j=2$.

\begin{figure}[h]
\centering
\includegraphics[width=0.5\textwidth]{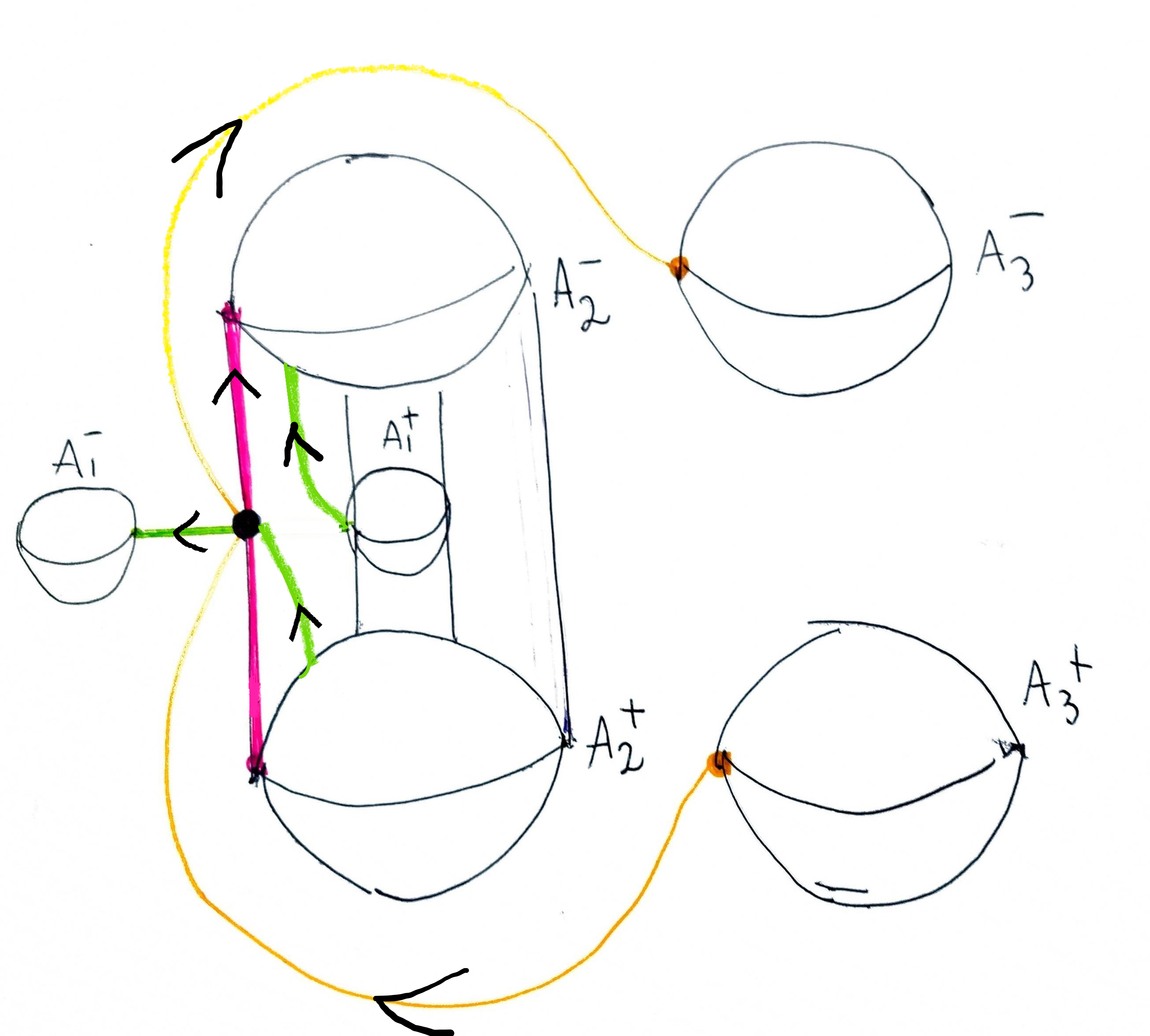}
\caption{The image of $a_1$ under $F_{1,2}$ is depicted in green, and it is homotopic to $a_1a_2$.}
\end{figure}
\FloatBarrier
Thus, $[F_{1,2}(a_1)]=[a_1a_2]$, and since $F_{1,2}$ fixes the homotopy classes of $a_2$ and $a_3$, then $F$ realizes $R_{1,2}$.

To see that $G_j$ actually realizes $I_j$, consider what $G_j$ does to $a_j$, as depicted in figure 4 for the case $j=1$. Notice that $G_j$ fixes the subpath of $a_j$ that is in $P_j\cap \{(x,y,z)|1/9\le x^2+y^2+z^2\le 1\}$.
\begin{figure}[hbt!]
\centering
\includegraphics[width=0.5\textwidth]{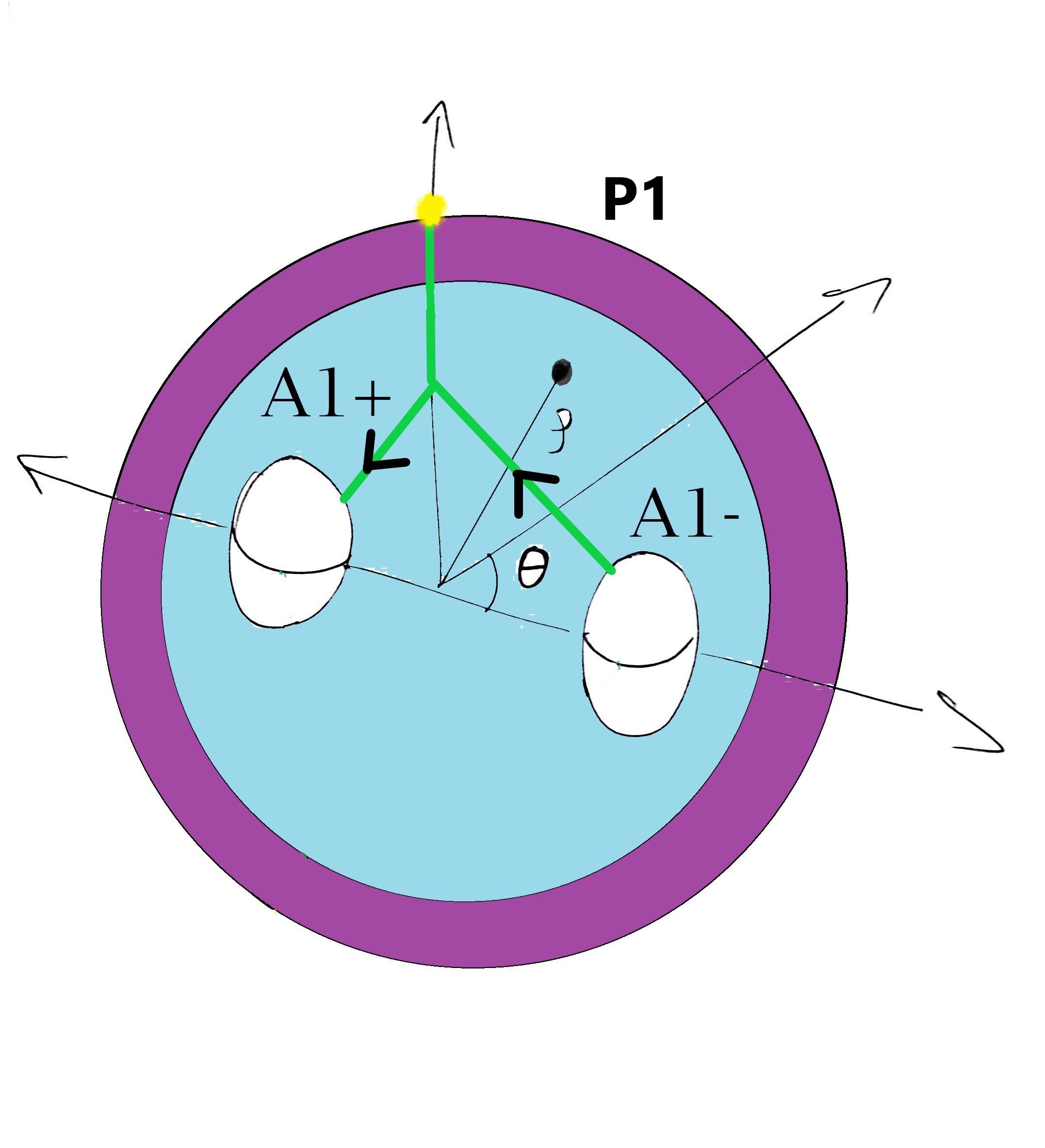}
\caption{The image of $a_1$ under $G_1$ is depicted in green, and it is homotopic to $a_1^{-1}$}.
\end{figure}
\FloatBarrier

Hence, $s$ defined on the Nielsen generators by $s([R_{i,j}])=[F_{i,j}]$, and $s([I_j])=G_j$ will be a section of $\rho$, provided that $\mathfrak{T}([F_{i,j}])=0$, and $\mathfrak{T}([G_{j}])=0$.

\section{Calculation of $\mathfrak{T}([F_{i,j}])$}

Denote  $N_{i,j}$  by $N$, $F_{i,j}$ by $F$, and  $f$ by $f_{i,j}$. Consider the universal cover $\widetilde{N}$ of $N$,  which is given by: 

$$\{(x,y,z)|x^{2}+y^{2}\le 1\}-\bigcup_{n\in\Z}C_n\subseteq \R^3,$$
where $C_n=\{(x,y,z)|x^{2}+y^{2}+(z-n)^{2}<1/9\}$.

Denote by $\pi$ the projection map $\pi:\widetilde{N}\xrightarrow{} N$ . $\pi$ is given by $\pi(x,y,z)=(x,y,e^{2\pi i z})$, and is a  local diffeomorphism. $f$ lifts to   a diffeomorphism $\widetilde{f}$ given by $\widetilde{f}(x,y,z)=(x,y,z+\psi\left(\sqrt{x^{2}+y^{2}}\right))$. 

We have that $\pi_1(\GL^{+}(3,\R),id)\cong \pi_1(SO(3),id) \cong \Z/2$ is generated by a loop $l:[0,1]\xrightarrow{} SO(3)$ which can be chosen to be  
$$l(t)=
\begin{bmatrix}
\cos(2\pi t) & -\sin(2\pi t) & 0\\
\sin(2\pi t) & \cos(2\pi t) & 0 \\
0 &0 & 1
\end{bmatrix},
$$

for $t\in [0,1]$.

For $M$ a closed oriented 3-manifold, let $TM$ be the tangent bundle of $M$ and define $Fr(TM)$ to be the principal $\GL_{3}^{+}(\R)$-bundle of oriented frames of $TM$. That means, $Fr(TM)_{x}$ is the space of linear isomorphism $T:\R^{3}\xrightarrow{} T_{x}M$. Fix a section $\sigma_0$ of $Fr(TM)$. We think of $\sigma_0$ as describing a preferred basis $\{\sigma_0(p)(e_1), \sigma_0(p)(e_2),\sigma_0(p)(e_3)\}$ of the tangent space at each point $p$.  Denote by $C(M,\GL^{+}(3,\R))$ the space of continuous functions from $M$ to $\GL^{+}(3,\R)$. 

The derivative crossed homomorphism $$\mathcal{D}: \text{Diff}^{+}(M)\xrightarrow{} C(M,\GL^{+}(3,\R))$$ will be defined now.  Given a diffeomrphism $F$ of $M$, the derivative crossed homomorphism evaluated at $[F]$, $\mathcal{D}([F]):M\xrightarrow{}\GL^{+}(3,\R)$, gives for each $p$ a linear transformation $\mathcal{D}([F])(p)$ in $\GL^{+}(3,\R)$, defined as follows.  It is the unique linear transformation that makes the following diagram commute:

\[ \begin{tikzcd}[row sep=25,column sep = 40]
\R^{3} \arrow{r}{\sigma_0(p)} \arrow[swap]{d}{\mathcal{D} ([F])(p)} & T_pN \\
\R^{3} \arrow{r}{\sigma_0(F(p))}& T_{F(p)}N \arrow{u}[swap] {[DF^{-1}]_{F(p)}}
\end{tikzcd}
\]

 Thus, $\mathcal{D}([F])(p)$ is the inverse of the linear transformation that represents the change of basis transformation from the basis

$\{\sigma_0(p)(e_1), \sigma_0(p)(e_2),\sigma_0(p)(e_3)\}$ of $T_{p}N$ to the basis  

$\{DF^{-1}(\sigma_0(F(p))(e_1)), DF^{-1}(\sigma_0(F(p))(e_2)), DF^{-1}(\sigma_0(F(p))(e_3))\}$  of $T_pN$, as depicted in figure 5.

\begin{figure}[hbt!]
\centering
\includegraphics[width=0.5\textwidth]{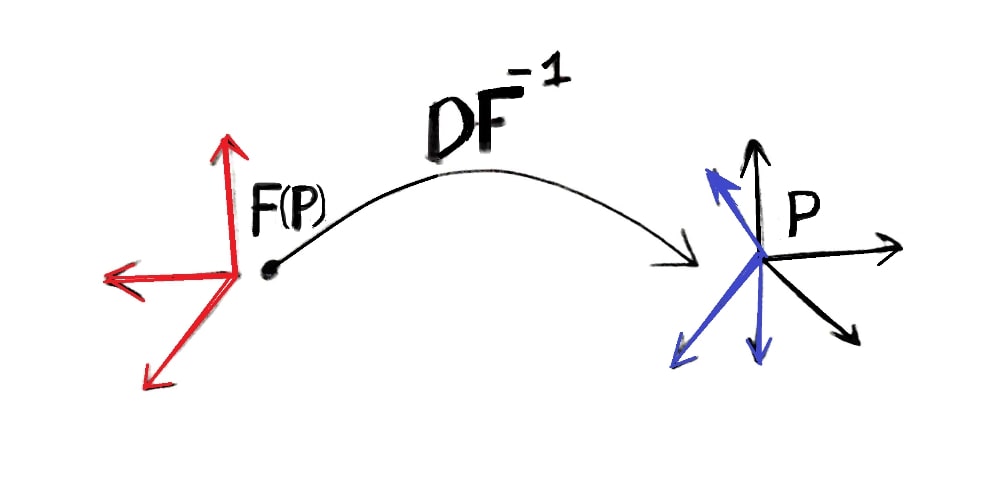}
\caption{The basis at $F(p)$ and $p$ that are determined by $\sigma_0$ are depicted in red and black respectively. The  basis in red  is sent to the basis in blue at $P$ by $DF^{-1}$. $\mathcal{D}([F])^{-1}(p)$ is the change of basis between the blue and black basis }.
\end{figure}

Thus, we get:
$$\mathcal{D}([F])^{-1}(p)=\sigma_0^{-1}(p)[DF^{-1}]_{F(p)}\sigma_0(F(p)).$$

In the particular case  of $M=M_n$, we study the derivative  crossed homomorphism of $F$ using a lift of it on the universal cover of $N$. This simplifies the computation of the derivative crossed homomorphism of $F$.
Let $q$ be in the interior of $\widetilde{N}$.  There is an isomorphism of vector spaces $b_q:\R^3\xrightarrow{} T_q\widetilde{N}\cong T_q\R^{3}$, defined by $b_q(e_1)=\dfrac{\partial}{\partial x}\Bigr|_{q}, b_q(e_2)=\dfrac{\partial}{\partial y}\Bigr|_{q}, \text{ and }, b_q(e_3)=\dfrac{\partial}{\partial z}\Bigr|_{q}$, where $\{e_1,e_2,e_3\}$ is the standard basis of $\R^{3}$.  
Then, define $\sigma(q)\in Fr(T\widetilde{N})$ by $\sigma(q)=b_q$.  Since $\pi$ is a local diffeomorphism, it induces an isomorphism of vector spaces $D\pi_q:T_q\widetilde{N}\xrightarrow{} T_{\pi(q)} N$. Let $p$ be in the interior of $N$. Select any $q$ in the interior of $\widetilde{N}$ such that $\pi(q)=p$. Then  $\sigma_0:\R^{3}\xrightarrow{} T_{p}N$ is defined by $\sigma_0(p):=D\pi_{q}\circ \sigma(q)$. 
We want to show that $\sigma_0$ doesn't depend on the lift $q$ of $p$. Let $q^{\prime}$ be another lift of $p$. Consider the Deck transformation $\Gamma$ of $\widetilde{N}$ that sends $q^{\prime}$ to $q$, and is given by $\Gamma(x,y,z)=(x,y,z+ k)$, for some $k\in\Z$. Then, $D\Gamma_{q^{\prime}}\left(\dfrac{\partial}{\partial x_i}\Bigr|_{q^{\prime}}\right)=\dfrac{\partial}{\partial x_i}\Bigr|_{q}$. Since $\Gamma$ is a Deck transformation of $\widetilde{N}$, it satisfies $\pi\circ \Gamma=\pi$. Then, $D\pi_{q}\circ D\Gamma_{q^{\prime}} =D\pi_{q^{\prime}}$. 
Hence:
$$\left[D\pi_{q^{\prime}}\circ \sigma(q^{\prime})\right](e_i)= 
\left[D\pi_{q^{\prime}}\right](\sigma(q^{\prime})(e_i))=
\left[D\pi_q\circ D \Gamma_{q^{\prime}}\right]\left(\dfrac{\partial}{\partial x_i}\Bigr|_{q^{\prime}}\right)= D\pi_{q} \left(\dfrac{\partial}{\partial x_i}\Bigr|_{q} \right)= $$

$$\left[D\pi_{q}\circ \sigma(q)\right](e_i).$$

Thus, $D\pi_{q^{\prime}}\circ \sigma(q^{\prime})=D\pi_{q}\circ \sigma(q)$, so $\sigma_0(p)$ doesn't depend on the lift $q$ of $p$. Thus, $\sigma_0$ is in fact a smooth section $\sigma_0:Int(N)\xrightarrow{} Fr(T(Int(N)))$ of the frame bundle of $Int(N)$. 

$$ $$

\begin{lema}
    
Let $p\in Int(N)$, and $q\in \widetilde{N}$ with $\pi(q)=p$. 
Then, ${\mathcal{D}([F])^{-1}_{ki}}(p)= \dfrac{\partial \widetilde{f_k}^{-1}}{\partial x_i}\Bigr|_{\widetilde{f}(q)}$.

\begin{proof}

Since $\pi\circ \widetilde{f}=f\circ \pi$, then by the chain rule we get:

$$ D\pi_{\widetilde{f}(q)}\circ D\widetilde{f}_{q} = Df_{\pi(q)}\circ D\pi_{q}.$$

Since $D\pi_q:T_{q}\widetilde{N}\xrightarrow{} T_pN$ is a linear isomorphism for each $q$, then:

$$D\pi_{\widetilde{f}(q)}\circ D\widetilde{f}_{q}\circ \left[D\pi_{q}\right]^{-1} = Df_{p},$$

and thus:

$$D\pi_{q}\circ D\widetilde{f}^{-1}_{\widetilde{f}(q)}\circ  \left[D\pi_{\widetilde{f}(q)}\right]^{-1} = Df^{-1}_{f(p)}.$$

Thus:

$\sigma_0^{-1}(p)Df^{-1}_{f(p)}\sigma_0(f(p))= \sigma_0^{-1}(p)\left[D\pi_{q}\circ D\widetilde{f}^{-1}_{\widetilde{f}(q)}\circ  \left[D\pi_{\widetilde{f}(q)}\right]^{-1}\right]\sigma_0(f(p))=\sigma^{-1}(q)D\widetilde{f}^{-1}_{\widetilde{f}(q)}\sigma(\widetilde{f}(q))$.

Therefore:

$$\mathcal{D}([F])^{-1}(p) =\sigma^{-1}(q)D\widetilde{f}^{-1}_{\widetilde{f}(q)}\sigma(\widetilde{f}(q)).$$

For $p$ in the interior of $N$, and $q$ with $\pi(q)=p$, evaluation of $e_i$ produces:

$$\sigma^{-1}(q)D\widetilde{f}^{-1}_{\widetilde{f}(q)}\sigma(\widetilde{f}(q))(e_i)=\sigma^{-1}(q)D\widetilde{f}^{-1}_{\widetilde{f}(q)}(\dfrac{\partial }{\partial x_i}\Bigr|_{\widetilde{f}(q)})=\sigma^{-1}(q) \left(\dfrac{\partial \widetilde{f}^{-1}}{\partial x_i}\Bigr|_{\widetilde{f}(q)}\right)=$$

$$\sigma^{-1}(q)\left( \sum_{k=1}^{3}\left(\dfrac{\partial \widetilde{f_k}^{-1}}{\partial x_i}\Bigr|_{\widetilde{f}(q)}\right)\dfrac{\partial  }{\partial x_k}\Bigr|_{q}\right)=\sum_{k=1}^{3}\left(\dfrac{\partial \widetilde{f_k}^{-1}}{\partial x_i}\Bigr|_{\widetilde{f}(q)}\right) \sigma^{-1}(q)\left( \dfrac{\partial  }{\partial x_k}\Bigr|_{q}\right)=$$

$$ \sum_{k=1}^{3}\left(\dfrac{\partial \widetilde{f_k}^{-1}}{\partial x_i}\Bigr|_{\widetilde{f}(q)}\right) e_k.$$

 Thus:

$$ {\mathcal{D}([F])^{-1}_{ki}}(p)= \dfrac{\partial \widetilde{f_k}^{-1}}{\partial x_i}\Bigr|_{\widetilde{f}(q)}$$
\end{proof}
\end{lema}

This lemma permits the calculation of $\mathcal{D}([F])$ in terms of $\widetilde{f}$ as we mentioned. 

Now, we define the twisting crossed homomorphism 

$$\mathfrak{T}:Mod(M)\xrightarrow{}Hom(\pi_1(M,x_0), \pi_1(\GL^{+}(3,\R), id))=H^{1}(M; \Z/2).$$

The twisting crossed homomorphism evaluated at $[F]$, $\mathfrak{T}([F]):\pi_1(M,x_0)\xrightarrow{}\pi_1(\GL^{+}(3,\R), id)$,   is the homomorphism that sends $[\gamma]\in \pi_1(M,x_0)$ to the class $\mathfrak{T}(F)[\gamma]=[\mathcal{D}([F])(\gamma)]$, where  $\mathcal{D}([F])(\gamma)$ is the loop:

$$[0,1]\xrightarrow{\gamma}M\xrightarrow{\mathcal{D}([F])}\GL^{+}(3,\R). $$

$\mathfrak{T}([F])$ studies the induced map of $\mathcal{D}([F])$ on the fundamental groups.

Because the derivative of $F$ is the identity on $a_k$ for $k\neq i$, we get that
 $\mathfrak{T}([F])[a_k]$ is trivial   in $\pi_{1}(\GL^{+}(3,\R),id)$ for every $k\neq i$.
 
Choose  $\gamma\in[a_i]$, and one of its lifts $\widetilde{\gamma}$,  such that $\widetilde{\gamma}$ intersects $\widetilde{N}$ as $(s,0, 0)$, for $s\in[0,1]$. For $t\in[0,1]$ satisfying that  $\gamma(t)\notin Int(N)$, we have that  the derivative of $F$ is trivial at $\gamma(t)$, thus $\mathcal{D}([F])(\gamma(t))$ is the trivial matrix for such  $t$. Hence, we are only interested  in the case $\gamma(t)\in Int(N)$, and in this case $F=f$.

Notice that $\widetilde{f}^{-1}(x,y,z)=(x,y,z-\psi(\sqrt{x^{2}+y^{2}}))$. Then we obtain :
$${\mathcal{D}([F])^{-1}}(\gamma(t))=
\begin{bmatrix}
1 &  0 & 0 \\
0 & 1 & 0\\
\dfrac{-x}{\sqrt{x^{2}+y^{2}}}\dfrac{d\psi}{d r}(\sqrt{x^{2}+y^{2}}) & \dfrac{-y}{\sqrt{x^{2}+y^{2}}}\dfrac{d\psi}{dr}(\sqrt{x^{2}+y^{2}}) & 1
\end{bmatrix}
(s,0, \psi(s))=$$

$$
\begin{bmatrix}
1 &  0 & 0 \\
0 & 1 & 0\\
-\dfrac{d\psi(s)}{d r} & 0 & 1
\end{bmatrix}.
$$

Hence,

$${\mathcal{D}([F])^{-1}}(\gamma(t))=
\begin{bmatrix}
1 &  0 & 0 \\
0 & 1 & 0\\
-\dfrac{d\psi(s)}{d r} & 0 & 1
\end{bmatrix}.
$$

So,

$${\mathcal{D}([F])}(\gamma(t))=
\begin{bmatrix}
1 &  0 & 0 \\
0 & 1 & 0\\
\dfrac{d\psi(s)}{dr} & 0 & 1
\end{bmatrix}.
$$

For $(s,0,0)$, $0\le s\le 1$.

We have an homotopy from the trivial path to  this path,

$$ H:[0,1]\times [0,1]\xrightarrow{} \GL^{+}(3,\R), $$

given by: $$H(s,t)=\begin{bmatrix}
1 &  0 & 0 \\
0 & 1 & 0 \\
 t\dfrac{d\psi(s)}{d r} & 0 & 1
\end{bmatrix}.$$

Therefore, $[\mathcal{D}([F])(\gamma)]=1$ for every $t\in Dom(\gamma)$.

Hence, $\mathfrak{T}([F])[a_i]$ is trivial in $\pi_{1}(\GL^{+}(3,\R),id)$.

Therefore, the twisting crossed homomorphism $\mathfrak{T}$ evaluated  at $F$, $\mathfrak{T}([F])$, is trivial.

\section{Calculation of $\mathfrak{T}([G_j])$}
Now, we analyse the case of $G_j$ . Denote $G_j$ by $G$, and $P_j$ by $P$. Given $\gamma\in [a_j]$, define $\gamma_1(t):=\gamma(t/3)$ , $\gamma_2(t):=\gamma(1/3+t/3)$ and $\gamma_3(t):=\gamma(2/3+t/3)$ for $0\le t\le 1$.  Then, $\gamma=\gamma_1*\gamma_2*\gamma_3$. Homotope $\gamma$ such that $\gamma_2\subseteq \{(x,y,z)|x^2+y^2+z^2\le 1/9\}$, and that $\gamma_1(t)=(0,0,s(t))$ and $\gamma_3(t)=(0,0,s(1-t))$, for some smooth function $s:[0,1]\xrightarrow{} [0,1]$. 
For $(x,y,z)\in P_j,$  $G$ is the diffeomorphism that sends $(x,y,z)$ to  

$$  \begin{bmatrix}
\cos(\psi(\sqrt{x^{2}+y^{2}+z^{2}})\pi) &  -\sin(\psi(\sqrt{x^{2}+y^{2}+z^{2}})\pi) & 0 \\
\sin(\psi(\sqrt{x^{2}+y^{2}+z^{2}})\pi) & \cos(\psi(\sqrt{x^{2}+y^{2}+z^{2}})\pi) & 0\\
0 & 0 & 1
\end{bmatrix}  \begin{pmatrix}
x \\
y \\
z
\end{pmatrix}. $$  

Evaluating  $G$ on $\gamma_1$  we get:

$$ G(\gamma_1(t))=\begin{bmatrix}
\cos(\psi(s(t))\pi) &  -\sin(\psi(s(t))\pi) & 0 \\
\sin(\psi(s(t))\pi) & \cos(\psi(s(t))\pi) & 0\\
0 & 0 & 1
\end{bmatrix}  \begin{pmatrix}
0 \\
0 \\
s(t)
\end{pmatrix}=(0,0,s(t)), $$ 

for $0\le t\le 1$.

  Evaluating $G$ on $\gamma_2$  we get:
$$ G(\gamma_2(t))= \begin{bmatrix}
-1 & 0 & 0 \\
0 & -1 & 0 \\
0 & 0 & 1
\end{bmatrix} \gamma_2(t)=\overline{\gamma_2(t)}, $$ 

for $0\le t\le 1$.

Evaluating $G$ on $\gamma_3$  we get:

$$  G(\gamma_3(t))=\begin{bmatrix}
\cos(\psi(s(1-t))\pi) &  -\sin(\psi(s(1-t))\pi) & 0 \\
\sin(\psi(s(1-t))\pi) & \cos(\psi(s(1-t))\pi) & 0\\
0 & 0 & 1
\end{bmatrix}  \begin{pmatrix}
0 \\
0 \\
s(1-t)
\end{pmatrix}=(0,0,s(1-t)), $$ 

for $0\le t\le 1$.

 Denote ${\mathcal{D}([G])^{-1}}(\gamma_1(t))$ by $w_1(t)$ and ${\mathcal{D}([G])^{-1}}(\gamma_3(t))$ by $w_3(t)$,  for $0\le t \le 1$.  The fact that $G(\gamma_1(t))=\overline{G(\gamma_3(t))}$ , implies that $w_1(t)=\overline{w_3(t)}$ for $0\le t\le 1.$ 

On the other hand, ${\mathcal{D}([G])^{-1}}(\gamma_2(t))$ is  the constant path:

$$ {\mathcal{D}([G])^{-1}}(\gamma_2(t))= \begin{bmatrix}
-1 &  0 & 0 \\
0 & -1 & 0\\
0 & 0 & 1
\end{bmatrix}=:w_2(t),$$   for $0\le t\le 1$

Therefore, ${\mathcal{D}([G])^{-1}}(\gamma)=w_1*w_2*\overline{w_1}$, which implies that $[{\mathcal{D}([G])^{-1}}(\gamma)]$ is trivial for $t\in Dom(\gamma)$, because $w_2(t)$ is a constant path. Therefore, the twisting crossed homomorphism $\mathfrak{T}$ evaluated  at $G$, $\mathfrak{T}([G])$, is trivial.

 \section{References}
 $$ $$
 
[1] F. Laudenbach, Sur les 2-sphres d’une variete de dimension 3, Ann. of Math. (2) 97 (1973), 57–81.

[2] F. Laudenbach, Topologie de la dimension trois: homotopie et isotopie, Societe Mathematique de France, Paris, 1974

[3] Brendle, T., Broaddus, N., $\&$ Putman, A. (2023). The mapping class group of connect sums of $S^2\times S^1$. Transactions of the American Mathematical Society, 376(04), 2557-2572.

[4]  R. C. Lyndon and P. E. Schupp, Combinatorial Group Theory. Berlin, Heidelberg, New York: Springer, 1977. 

\end{document}